\documentclass{amsart}
\usepackage{amsmath,amssymb,amsbsy,amsfonts,amsthm,latexsym,amsopn,amstext,
                                               amsxtra,euscript,amscd}
\begin{document}

\newfont{\teneufm}{eufm10}
\newfont{\seveneufm}{eufm7}
\newfont{\fiveeufm}{eufm5}
%
%
\newfam\eufmfam
                \textfont\eufmfam=\teneufm \scriptfont\eufmfam=\seveneufm
                \scriptscriptfont\eufmfam=\fiveeufm
%
%
\def\frak#1{{\fam\eufmfam\relax#1}}
%


\def\bbbr{{\rm I\!R}} 
\def\bbbm{{\rm I\!M}}
\def\bbbn{{\rm I\!N}} 
\def\bbbf{{\rm I\!F}}
\def\bbbh{{\rm I\!H}}
\def\bbbk{{\rm I\!K}}
\def\bbbp{{\rm I\!P}}
\def\bbbone{{\mathchoice {\rm 1\mskip-4mu l} {\rm 1\mskip-4mu l}
{\rm 1\mskip-4.5mu l} {\rm 1\mskip-5mu l}}}
\def\bbbc{{\mathchoice {\setbox0=\hbox{$\displaystyle\rm C$}\hbox{\hbox
to0pt{\kern0.4\wd0\vrule height0.9\ht0\hss}\box0}}
{\setbox0=\hbox{$\textstyle\rm C$}\hbox{\hbox
to0pt{\kern0.4\wd0\vrule height0.9\ht0\hss}\box0}}
{\setbox0=\hbox{$\scriptstyle\rm C$}\hbox{\hbox
to0pt{\kern0.4\wd0\vrule height0.9\ht0\hss}\box0}}
{\setbox0=\hbox{$\scriptscriptstyle\rm C$}\hbox{\hbox
to0pt{\kern0.4\wd0\vrule height0.9\ht0\hss}\box0}}}}
\def\bbbq{{\mathchoice {\setbox0=\hbox{$\displaystyle\rm
Q$}\hbox{\raise 0.15\ht0\hbox to0pt{\kern0.4\wd0\vrule
height0.8\ht0\hss}\box0}} {\setbox0=\hbox{$\textstyle\rm
Q$}\hbox{\raise 0.15\ht0\hbox to0pt{\kern0.4\wd0\vrule
height0.8\ht0\hss}\box0}} {\setbox0=\hbox{$\scriptstyle\rm
Q$}\hbox{\raise 0.15\ht0\hbox to0pt{\kern0.4\wd0\vrule
height0.7\ht0\hss}\box0}} {\setbox0=\hbox{$\scriptscriptstyle\rm
Q$}\hbox{\raise 0.15\ht0\hbox to0pt{\kern0.4\wd0\vrule
height0.7\ht0\hss}\box0}}}}
\def\bbbt{{\mathchoice {\setbox0=\hbox{$\displaystyle\rm
T$}\hbox{\hbox to0pt{\kern0.3\wd0\vrule height0.9\ht0\hss}\box0}}
{\setbox0=\hbox{$\textstyle\rm T$}\hbox{\hbox
to0pt{\kern0.3\wd0\vrule height0.9\ht0\hss}\box0}}
{\setbox0=\hbox{$\scriptstyle\rm T$}\hbox{\hbox
to0pt{\kern0.3\wd0\vrule height0.9\ht0\hss}\box0}}
{\setbox0=\hbox{$\scriptscriptstyle\rm T$}\hbox{\hbox
to0pt{\kern0.3\wd0\vrule height0.9\ht0\hss}\box0}}}}
\def\bbbs{{\mathchoice
{\setbox0=\hbox{$\displaystyle     \rm S$}\hbox{\raise0.5\ht0\hbox
to0pt{\kern0.35\wd0\vrule height0.45\ht0\hss}\hbox
to0pt{\kern0.55\wd0\vrule height0.5\ht0\hss}\box0}}
{\setbox0=\hbox{$\textstyle        \rm S$}\hbox{\raise0.5\ht0\hbox
to0pt{\kern0.35\wd0\vrule height0.45\ht0\hss}\hbox
to0pt{\kern0.55\wd0\vrule height0.5\ht0\hss}\box0}}
{\setbox0=\hbox{$\scriptstyle      \rm S$}\hbox{\raise0.5\ht0\hbox
to0pt{\kern0.35\wd0\vrule height0.45\ht0\hss}\raise0.05\ht0\hbox
to0pt{\kern0.5\wd0\vrule height0.45\ht0\hss}\box0}}
{\setbox0=\hbox{$\scriptscriptstyle\rm S$}\hbox{\raise0.5\ht0\hbox
to0pt{\kern0.4\wd0\vrule height0.45\ht0\hss}\raise0.05\ht0\hbox
to0pt{\kern0.55\wd0\vrule height0.45\ht0\hss}\box0}}}}
\def\bbbz{{\mathchoice {\hbox{$\sf\textstyle Z\kern-0.4em Z$}}
{\hbox{$\sf\textstyle Z\kern-0.4em Z$}} {\hbox{$\sf\scriptstyle
Z\kern-0.3em Z$}} {\hbox{$\sf\scriptscriptstyle Z\kern-0.2em
Z$}}}}
\def\ts{\thinspace}

\newtheorem{theorem}{Theorem}
\newtheorem{lemma}[theorem]{Lemma}
\newtheorem{claim}[theorem]{Claim}
\newtheorem{cor}[theorem]{Corollary}
\newtheorem{prop}[theorem]{Proposition}
\newtheorem{definition}{Definition}
\newtheorem{question}[theorem]{Open Question}

\def\squareforqed{\hbox{\rlap{$\sqcap$}$\sqcup$}}
\def\qed{\ifmmode\squareforqed\else{\unskip\nobreak\hfil
\penalty50\hskip1em\null\nobreak\hfil\squareforqed
\parfillskip=0pt\finalhyphendemerits=0\endgraf}\fi}

\def\cA{{\mathcal A}}
\def\cB{{\mathcal B}}
\def\cC{{\mathcal C}}
\def\cD{{\mathcal D}}
\def\cE{{\mathcal E}}
\def\cF{{\mathcal F}}
\def\cG{{\mathcal G}}
\def\cH{{\mathcal H}}
\def\cI{{\mathcal I}}
\def\cJ{{\mathcal J}}
\def\cK{{\mathcal K}}
\def\cL{{\mathcal L}}
\def\cM{{\mathcal M}}
\def\cN{{\mathcal N}}
\def\cO{{\mathcal O}}
\def\cP{{\mathcal P}}
\def\cQ{{\mathcal Q}}
\def\cR{{\mathcal R}}
\def\cS{{\mathcal S}}
\def\cT{{\mathcal T}}
\def\cU{{\mathcal U}}
\def\cV{{\mathcal V}}
\def\cW{{\mathcal W}}
\def\cX{{\mathcal X}}
\def\cY{{\mathcal Y}}
\def\cZ{{\mathcal Z}}

\newcommand{\comm}[1]{\marginpar{%
\vskip-\baselineskip 
\raggedright\footnotesize
\itshape\hrule\smallskip#1\par\smallskip\hrule}}




\newcommand{\ignore}[1]{}

\def\vec#1{\mathbf{#1}}

\def\e{\mathbf{e}}



\def\GL{\mathrm{GL}}

\hyphenation{re-pub-lished}

\def\rank{{\mathrm{rk}\,}}
\def\dd{{\mathrm{dyndeg}\,}}

\def\A{\mathbb{A}}
\def\B{\mathbf{B}}
\def \C{\mathbb{C}}
\def \F{\mathbb{F}}
\def \K{\mathbb{K}}
\def \Z{\mathbb{Z}}
\def \P{\mathbb{P}}
\def \R{\mathbb{R}}
\def \Q{\mathbb{Q}}
\def \N{\mathbb{N}}
\def \Z{\mathbb{Z}}

\def \nd{{\, | \hspace{-1.5 mm}/\,}}

\def\Zn{\Z_n}

\def\Fp{\F_p}
\def\Fq{\F_q}
\def \fp{\Fp^*}
\def\\{\cr}
\def\({\left(}
\def\){\right)}
\def\fl#1{\left\lfloor#1\right\rfloor}
\def\rf#1{\left\lceil#1\right\rceil}

\setlength{\textheight}{43pc}
\setlength{\textwidth}{28pc}

\title[Polynomial Dynamics and Nonlinear Pseudorandom Numbers]
{On the Degree Growth in Some Polynomial Dynamical
Systems and Nonlinear Pseudorandom Number Generators}

\author{Alina~Ostafe}
\address{Institut f\"ur Mathematik, Universit\"at Z\"urich,
Winterthurerstrasse 190 CH-8057, Z\"urich, Switzerland}
\email{alina.ostafe@math.uzh.ch}

\author{Igor E.~Shparlinski} 
\address{Department of Computing, Macquarie University\\ NSW 2109, Australia}
\email{igor@ics.mq.edu.au}

\begin{abstract}  In this paper we study a class 
of dynamical systems generated by iterations of multivariate polynomials  
and estimate the degree
growth of these iterations. We use these estimates to 
bound exponential sums along the 
orbits of these dynamical systems and show that they admit 
much stronger estimates
than in the general case and thus can be of use for pseudorandom number 
generation. 
\end{abstract}

\maketitle


\section{Introduction}

Given a system of $r$ polynomials $\cF = \{f_0, \ldots ,f_{r-1}\}$ 
in $r$ variables over a
ring $\cR$ one can naturally define a  dynamical system
generated by its iterations:
$$
f_i^{(0)}=f_i, \qquad f_i^{(k)}= f_i^{(k-1)}(f_0, \ldots ,f_{r-1}), \qquad k=0,1, \ldots\ ,  
$$
for each $i = 0, \ldots, r-1$, 
see~\cite{Arn1,Arn2,Arn3,Chang,CohHac,CJLS,EvWa,FomZel,Jon,Marc,Shp,Silv1,Silv2}
and references therein 
for various aspects of such dynamical systems. 
It is also natural to consider the orbits obtained by such 
iterations evaluated at a certain initial value 
$\(u_{k,0}, \ldots, u_{k,r-1}\)$.

In the special case of one linear univariate polynomial over 
a residue ring or a finite field such iterations, 
known as linear congruential generators, have been successully 
used for decades in the theory of Quasi Monte Carlo methods, see~\cite{Nied1,Nied2}.

Unfortunately, in cryptographic settings, such linear generators 
have been successfully attacked~\cite{ContShp,FHKLS,JS,Kraw,Lag} and thus deemed unusable 
for cryptograpic purposes. 
It should be noted that nonlinear generators have also been attacked~\cite{BGGS1,BGGS2,GoGuIb,GuIb},
but the attacks are much weaker and do not rule out their 
use for cryptographic purposes (provided reasonable precausions are 
made). 
Although linear congruential generators have been used quite sucessfully
for  Quasi-Monte Carlo methods, their linear structure shows in 
these applications too and often limits their applicability,
see~\cite{Nied1,Nied2}. 

 Motivated by these potential applications, 
the statistical uniformity of the distribution 
(measured by the discrepancy) of one and  multidimensional 
nonlinear polynomial generators 
have been studied in~\cite{GNS,GG,NiSh1,NiSh4,NiWi,TopWin}. However, all
previously known results are nontrivial only for those 
polynomial generators
that produce sequences of extremely large period, which could
be hard to achieve in practice. The reason
behind this is that typically the degree of iterated polynomial
systems grows exponentially, and that in all previous results
the saving over the trivial bound has been logarithmic.
Furthermore, it is easy to see that in the one dimensional
case (that is, for $r=1$) the exponential growth of 
the degree of iterations of a nonlinear polynomial is
unavoidable. One also expects the same behaviour 
in the mulitidimensional case for ``random'' polynomials
$f_0, \ldots ,f_{r-1}$. However, for some specially selected
polynomials $f_0, \ldots ,f_{r-1}$ the degree may grow 
significantly slower.

Indeed, here we describe a rather wide class of 
polynomial systems with polynomial growth of the degree
of their iterations. As a result we obtain much better estimates
of exponential sums, and thus of discrepancy, for vectors 
generated by these iterations, with a saving over 
the trivial bound  being a power of $p$. 
Our construction resembles that of {\it triangular maps\/} of~\cite{Marc} but behaves 
quite differently; for example, triangular maps in~\cite{Marc}
have the fastest possible degree growth. 
 
We remark that in the case of the so-called {\it inversive generator\/}
rather strong estimates are also available~\cite{NiSh2,NiSh3}, but this generator 
involves a modular inversion at each step which is a computationally 
expensive operation.  Another alternative where stronger than general
bounds are known is the power generator which essentially 
consists of iterating a monomial map $X\to X^e$, 
see~\cite{ChouShp,FHS,FPS,FrSh,MartPom,VaSha} and especially 
the recent result of J.~Bourgain~\cite{Bour} on 
the joint distribution of consecutive terms of this generator.
Similar results also hold for pseudorandom number  
generators obtained by iterating 
{\it Dickson polynomials\/}~\cite{GoGuSh} 
and {\it Redei functions\/}~\cite{GuWin}. 

Finally, we note that we also hope that our results may be of use 
for some applications in polynomial dynamical systems.

\section{Polynomial Dynamical System with Slow Degree Growth}

\subsection{Construction}

Let  $\F$ be an arbitrary field of characteristic $p>m$ (or of zero characteristic) 
and let $\cF = \{f_0, \ldots ,f_m\}$ be a system 
of $m+1$ polynomials in $\F[X_0, \ldots ,X_m]$ defined in the following way:
\begin{equation}
\label{eq:Polys}
\begin{split}
f_0(X_0, \ldots ,X_m)& = X_0g_0(X_1,\ldots,X_m)+h_0(X_1,\ldots,X_m),\\
f_1(X_0, \ldots ,X_m)&=X_1g_1(X_2,\ldots,X_m)+h_1(X_2,\ldots,X_m),\\
  &\ldots  \\
f_{m-1}(X_0, \ldots ,X_m)&=X_{m-1}g_{m-1}(X_m)+h_{m-1}(X_m),\\
f_m(X_0, \ldots ,X_m) &=aX_m+b,
\end{split}
\end{equation}
where $$
a,b\in\F,\quad a \ne 0, \quad g_i,h_i\in\F[X_{i+1},\ldots,X_m],\quad i=0,\ldots,m-1.
$$ We also impose the condition that each polynomial $g_i$ has {\it unique leading monomial\/} $X_{i+1}^{s_{i,i+1}}\ldots X_m^{s_{i,m}}$, that is,
\begin{equation}
\label{eq:Cond1}
g_i(X_{i+1},\ldots,X_m) = X_{i+1}^{s_{i,i+1}}\ldots X_m^{s_{i,m}} + 
\widetilde{g_i}(X_{i+1},\ldots,X_m), 
\end{equation}
where 
\begin{equation}
\label{eq:Cond2}
\deg \widetilde{g_i} < \deg g_i = s_{i,i+1}+\ldots+s_{i,m}
\end{equation}
and 
\begin{equation}
\label{eq:Cond3}
\deg h_i\le \deg g_i
\end{equation}
for $i=0,\ldots,m-1$.

For each $i=0, \ldots ,m$ we define the $k$-th iteration of the polynomials $f_i$ by the recurrence relation
\begin{equation}
\label{eq:PolyIter}
f_i^{(0)}=f_i, \qquad f_i^{(k)}= f_i^{(k-1)}(f_0, \ldots ,f_m), \qquad k=0,1, \ldots\,  .
\end{equation}

\subsection{Degree Growth} 

We  denote by $d_{k,i}$ the degree of the polynomial $f_i^{(k)}$, $i=0,\ldots,m$. We  also consider the vector of degrees of the iterations
$$\vec{d}_k=(d_{k,0} \ldots ,d_{k,m}),
$$ 
and the upper triangular matrix 
$$S=   \(\begin{matrix}        1 & s_{0,1} &  s_{0,2} & \ldots & s_{0,m} \\
      0 & 1 & s_{1,2} & \ldots & s_{1,m} \\
      &&\ldots &&\\
      0 & 0 & 0 & \ldots & 1
   \end{matrix}\)
 $$
 given by the exponents of the leading monomials in $f_i$, $i=0,\ldots,m$.
 We observe that under iterations we have 
 \begin{eqnarray*}
 f_i^{(k)}&=&f_i^{(k-1)}g_i(f_{i+1}^{(k-1)},\ldots,f_m^{(k-1)})+h_i(f_{i+1}^{(k-1)},\ldots,f_m^{(k-1)}),\\
 & &  \qquad \qquad \qquad \qquad \qquad \qquad \qquad \qquad \quad i=0,\ldots,m-1,\\
 f_m^{(k)}&=&af_m^{(k-1)}+b, 
 \end{eqnarray*}
 and using the conditions on the degrees of the polynomials $g_i$ and $h_i$ we get
 \begin{eqnarray*}
d_{k,i}&=&d_{k-1,i}+s_{i,i+1}d_{k-1,i+1}+\ldots+s_{i,m}d_{k-1,m},\qquad i=0,\ldots,m-1,\\
d_{k,m}&=&1.
\end{eqnarray*}
Using the above notations, the degrees of the iterations satisfy the relation
$$
\vec{d}_k=S \vec{d}_{k-1}, \qquad k\ge 0 \qquad \text{and} 
\qquad \vec{d}_{-1}=(1,\ldots,1)^t
$$
which is equivalent to writing 
\begin{equation}
\label{eq:dk explicit}
\vec{d}_k=S^{k+1}(1 \ldots ,1)^{t},\qquad k\ge0.
\end{equation}

We now show that the degrees of the iterations of $\cF$ grow polynomially.

\begin{lemma}
\label{lem:Poly Deg} Let $f_0, \ldots, f_m\in \F[X_0,\ldots,X_m]$ be  as in~\eqref{eq:Polys},
 satisfying the conditions~\eqref{eq:Cond1}, \eqref{eq:Cond2} and~\eqref{eq:Cond3}.
Then the degrees of the iterations of $\cF = \{f_0, \ldots ,f_m\}$ 
grow as follows
  \begin{eqnarray*} d_{k,i}&=&\frac{1}{(m-i)!}k^{m-i}s_{i,i+1}\ldots s_{m-1,m}+\psi_i(k),\qquad  i=0,\ldots, m-1,\\
d_{k,m}&=&1, 
\end{eqnarray*} 
where $\psi_i(T) \in \Q[T]$ is a polynomial of degree  $\deg \psi_i <m-i$.
 \end{lemma} 
 
\begin{proof} We  use induction on $m$. For $m=1$ one can easily see that we get
$$
d_{k,0}=ks_{0,1}+s_{0,1}+1\qquad \text{and}\qquad d_{k,1}=1.
$$
We assume the result true for $m$ indeterminates. 
Let $S$ be the matrix of exponents of the leading monomials in $\cF$ as above. We write $S$ in the following way
$$
S=   \(\begin{matrix} 
      R & \vec{s}\\
      0 & 1 
   \end{matrix}\),
$$
where $R$ is the matrix given by the exponents of the first $m$ indeterminates in the leading monomials of $f_i$, $i=0,\ldots,m-1$, and  $\vec{s}=(s_{0,m},\ldots,s_{m-1,m})$.
For a vector $\vec{v} \in \F^m$ we use $\vec{v}^t$ and $\vec{v}_i$ to denote the transpose
and the $i$th component  of $\vec{v}$, respectively.
We also denote by $\vec{e}$ the unit vector $\vec{e} = (1,\ldots,1) \in \F^m$.
Using these notations and recalling~\eqref{eq:dk explicit}, we obtain
$$
\vec{d}_k=S^{k+1}\vec{e}^t=   \(\begin{matrix} 
     R^{k+1} & (R^k+\ldots+R+I)\vec{s}^t \\
     0 & 1
   \end{matrix}\)\vec{e}^t.
$$
Componentwise, we have
\begin{eqnarray*}
d_{k,i}&=&\(R^{k+1}\vec{e}^t\)_i+\((R^k+\ldots+R+I)\vec{s}^t\)_i, \qquad 
i=0,\ldots, m-1,\\
d_{k,m}&=&1.
\end{eqnarray*}
It is easy to note that the maximal degree of the $k^{th}$-iteration of polynomials $f_i$ for any $i$ is given by the last position in each row of $S^{k+1}$. Using this remark and the induction hypothesis we get
$$
(R^j \vec{s}^t)_i
=\frac{1}{(m-1-i)!}j^{m-1-i}s_{i,i+1}\ldots s_{m-2,m-1}s_{m-1,m}+ \varphi_i(j), 
$$
for some polynomials $\varphi_i(Z) \in \Q[Z]$ of degree $\deg \varphi_i<m-1-i$.
Then
$$
\sum_{j=0}^k(R^j\vec{s}^t)_i=\frac{1}{(m-1-i)!}s_{i,i+1}\ldots s_{m-1,m}\sum_{j=0}^kj^{m-1-i}+\widetilde{\varphi}_i(k),
$$
for some polynomials $\widetilde{\varphi}_i(Z) \in \Q[Z]$ of degree $\deg \widetilde{\varphi}_i<m-i$.
As
$$
\sum_{j=0}^kj^{m-1-i}=\frac{1}{m-i}(B_{m-i}(k+1)-B_{m-i}(0)),
$$
where $B_{m-i}$ is  the Bernoulli polynomial of degree $m-i$ 
(which has the leading coefficient equal to $1$), we finally obtain the desired result. 
\end{proof}

\begin{cor}
\label{cor:Lin Indep}
Let $f_0, \ldots, f_m\in \F[X_0,\ldots,X_m]$ be  as in~\eqref{eq:Polys},
 satisfying the conditions~\eqref{eq:Cond1}, \eqref{eq:Cond2} and~\eqref{eq:Cond3}.
If $s_{0,1} \ldots s_{m-1,m} \ne 0$,  then for any integer $\nu \ge 1$
there is a constant $k_0$ depending only on the matrix $S$ and $\nu$ such that 
for any  integers $k_1,\ell_1, \ldots, k_\nu, \ell_\nu \ge k_0$ and 
any nonzero $\vec{a} = (a_0, \ldots, a_{m-1}) \in \F^m$, 
$$
F_{\vec{a},k_1,\ell_1, \ldots, k_\nu, \ell_\nu} =  \sum_{i=0}^{m-1} \!a_i
\sum_{j=1}^\nu \(f_i^{(k_j)}-f_i^{(\ell_j)}\), $$
is a nonconstant polynomial of degree 
$$
\deg F_{\vec{a}, k_1,\ell_1, \ldots, k_\nu, \ell_\nu}  = O(k^m),
$$ 
where 
$$
k = \max \{k_1,\ell_1, \ldots, k_\nu, \ell_\nu\}
$$
unless the components of the vectors 
 $$
 (k_1 \ldots, k_\nu)  \qquad \text{and}\qquad  (\ell_1 \ldots, \ell_\nu)
 $$
 are permutations of each other.
 \end{cor}
 
\begin{proof} 
Let  $i_0$ be the smallest integer with $a_{i_0} \ne 0$. Performing 
all trivial cancellations, without loss 
of generailty we can also assume that the vectors
$(k_1 \ldots, k_\nu)$ and $(\ell_1 \ldots, \ell_\nu)$ have no 
common elements. Thus the largest element amongst them $k$,
is unique. It is now clear from Lemma~\ref{lem:Poly Deg} that 
the leading term of $f_{i_0}^{(k)}$ is present in 
$F_{\vec{a},k_1,\ell_1, \ldots, k_\nu, \ell_\nu}$. 
\end{proof}
 
\section{Polynomial Pseudorandom Number Generators}

\subsection{Construction}

Let $\cF=\{f_0,\ldots,f_m\}$ be a list of $m+1$ polynomials in $\F_p[X_0,\ldots,X_m]$ defined as in section 2. 
We consider the sequence defined by a recurrence congruence
modulo a prime $p$ of the form
\begin{equation}
\label{eq:Gen}
u_{n+1,i}\equiv f_i(u_{n,0},\ldots,u_{n,m})\!\!\! \pmod p, \qquad n = 0,1,\ldots,
\end{equation}
with some {\it initial values\/}
$u_{0,0},\ldots,u_{0,m}$. 
We also assume that $0 \le u_{n,i} < p$, $i=0,\ldots,m$, $n=0, 1, \ldots$.
Using the following vector notation
$$
\vec{w_n}=(u_{n,0},\ldots,u_{n,m})
$$
and
$$
\cF=(f_0(X_0,\ldots,X_m),\ldots,f_m(X_0,\ldots,X_m)),
$$
we have the recurrence relation
$$
\vec{w_{n+1}}=\cF(\vec{w_{n}}).
$$
In particular, for any $n,k\ge 0$ and $i=0,\ldots,m$ we have
$$
u_{n+k,i}=f_i^{(k)}(u_{n,0},\ldots,u_{n,m})
$$
or 
$$
\vec{w_{n+k}}=\cF^{(k)}(\vec{w_n}).
$$
Clearly the sequence of vectors $\vec{w_{n}}$ is eventually periodic
with some period $T \le p^{m+1}$.
Without loss of generality we assume that it is 
$$
\vec{w_{n+T}} = \vec{w_{n}}, \qquad n =0,1, \ldots\,.
$$
In our construction of pseudorandom sequences,
we discard the last component in the vectors $\vec{w_{n}}$
and denote 
$$
\vec{u}_n=(u_{n,0},\ldots,u_{n,m-1})
$$
which we show to be rather uniformly distributed provided $T$ is large 
enough.

\subsection{Exponential Sums}
\label{sec:ExpSum}

We put
$$\e(z) = \exp(2 \pi i z/p).$$
Our second main tool is the Weil bound on exponential sums
(see~\cite[Chapter~5]{LN}) which we present
in the following slightly generalized form.

\begin{lemma}
\label{lem:Weil}
For any nonconstant  polynomial  $F \in
\F_p[X_{0},\ldots,X_{m}]$
of total degree $D$ we have the bound
$$
\left|\sum_{x_0, \ldots, x_m =1 }^p
\e\( F(x_{0},\ldots,x_{m})\) \right| < D p^{m+1/2}.
$$
\end{lemma}

We follow the  scheme previously introduced in~\cite{NiSh1,NiSh2}.
Furthermore, as it has been suggested in~\cite{NiWi,TopWin},
we  work with higher moments 
of the corresponding exponential sums. 
However  the polynomial growth of the
degree allows us a much more favorable choice of parameters and 
thus leads to a better estimate than in previous works.

Assume that  the sequence $\{\vec{u}_n\}$ generated by~\eqref{eq:Gen} is 
purely periodic
with an arbitrary period $T$.
For an integer vector $\vec{a} = (a_0, \ldots, a_{m-1}) \in \Z^m$ we introduce
the exponential sum
$$
S_{\vec{a}}(N) =  \sum_{n=0}^{N-1} \e\(\sum_{i=0}^{m-1} a_iu_{n,i}\).
$$

\begin{theorem}
\label{thm:ExpSum} 
Let  the sequence $\{\vec{u}_n\}$ be given by~\eqref{eq:Gen}, where
the family of $m+1$ polynomials 
$\cF=\{f_0,\ldots,f_m\} \in \F_p[X_{0},\ldots,X_{m}]$ 
of total degree $d \ge 2$ is  of the form~\eqref{eq:Polys}, 
satisfying the conditions~\eqref{eq:Cond1}, \eqref{eq:Cond2} and~\eqref{eq:Cond3}, 
 and such that $s_{0,1} \ldots s_{m-1,m} \ne 0$.
Assume that   $\{\vec{u}_n\}$ is purely
periodic with period $T$.  Then for any fixed integer $\nu\ge 1$, 
and any positive integer $N \le T$, 
the bound
$$
\max_{\gcd(a_0, \ldots, a_{m-1} ,p) = 1} \left| S_{\vec{a}}(N) \right| =
O\(p^{\alpha_{m,\nu}}N^{1-\beta_{m,\nu}}\)
$$
holds, where
$$
\alpha_{m,\nu}=\frac{2m^2 + 2m\nu + 2m +\nu}{4\nu(m+\nu)} \qquad \text{and} 
\qquad \beta_{m,\nu} = \frac{1}{2\nu}
$$
and the implied constant depends  only on $d$, $m$
and $\nu$.
\end{theorem}

\begin{proof} Select any $\vec{a} = (a_0, \ldots, a_{m-1}) \in \Z^m$  with
$\gcd(a_0, \ldots, a_{m-1} ,p) = 1$.
It is obvious that for any integer $k\ge 1$ we have
$$
\left|S_{\vec{a}}(N)   -  \sum_{n=0}^{N-1} \e\(\sum_{i=0}^{m-1} 
a_i u_{n+k,i}\)
\right| \le 2k.
$$
Let $k_0$ be the same as in Corollary~\ref{cor:Lin Indep}.
Therefore, for any integer $K \ge k_0$,
\begin{equation}
\label{eq:S and W}
(K-k_0+1) |S_{\vec{a}}(N)|  \le  W + K^2,
\end{equation}
where
$$
W = \left |\sum_{n=0}^{N-1}\sum_{k=k_0}^{K}
\e\(\sum_{i=0}^{m-1} a_i u_{n+k,i}\) \right|
 \le \sum_{n=0}^{N-1}\left | \sum_{k=k_0}^{K}
\e\(\sum_{i=0}^{m-1} a_i u_{n+k,i}\) \right| .
$$

As before, we define the sequence of polynomials 
$$f_i^{(k)}(X_{0},\ldots,X_{m}) \in
\F_p[X_{0},\ldots,X_{m}]$$
by~\eqref{eq:PolyIter}.  Then
\begin{eqnarray*}
W^{2\nu} & \le & N^{2\nu-1} \sum_{n=0}^{N-1}
\left | \sum_{k=k_0}^{K} \e\( \sum_{i=0}^{m-1} a_i
f_{i}^{(k)}\(\vec{u_n}\)\)
\right| ^{2\nu}\\
 & \le & N^{2\nu-1} \sum_{w_{0},\ldots, w_{m}\in \F_{p}}
\left | \sum_{k=k_0}^{K} \e\(\sum_{i=0}^{m-1} a_i
f_{i}^{(k)}\(w_{0},\ldots,w_{m}\) \) \right| ^{2\nu}\\
& = & N^{2\nu-1}
 \sum_{k_1,\ell_1, \ldots, k_\nu, \ell_\nu = k_0}^{K} 	\sum_{\vec{w}\in \F_{p}^{m+1}}
 	\e\!\(\sum_{i=0}^{m-1} \!a_i
	\sum_{j=1}^\nu 
 	\(f_i^{(k_j)}\(\vec{w}\) - f_i^{(\ell_j)}\(\vec{w}\)\)\).
\end{eqnarray*}
For  $O(K^\nu)$ vectors 
 $$
 (k_1 \ldots, k_\nu)  \qquad \text{and}\qquad  (\ell_1 \ldots, \ell_\nu)
 $$
which are permutations of each other,  we estimate  the inner sum 
trivially as $p^{m+1}$.

For the other $O(K^{2\nu})$ vectors, we combine  Corollary~\ref{cor:Lin Indep}
with Lemma~\ref{lem:Weil}  getting the
upper bound $K^mp^{m+1/2}$ for the inner sum
for at most $K^2$ sums. Hence, 
$$W^{2\nu} \le K^\nu N^{2\nu-1} p^{m+1} +  K^{m+2\nu} N^{2\nu-1} p^{m+1/2}.
$$
Inserting this bound in~\eqref{eq:S and W}, we derive
$$
S_{\vec{a}}(N) 
= O \(K^{-1/2}N^{1-1/2\nu} p^{(m+1)/2\nu} +  K^{m/2\nu} N^{1-1/2\nu} p^{(2m+1)/4\nu} +K\).
$$
Choosing
$$K = \rf{p^{1/2(m+\nu)}}$$
(and assuming that $p$ is large enough, so $K\ge k_0$),
after simple calculations  we obtain the desired result.
\end{proof}
 
Since 
$$
\lim_{\nu \to \infty} \alpha_{m,\nu}/\beta_{m,\nu} = m + 1/2
$$
we see from Theorem~\ref{thm:ExpSum} 
that for any fixed $\varepsilon> 0$ there is there $\delta>0$ 
such if 
$T\ge N \ge p^{m+1/2 + \varepsilon}$ 
then 
$$
\max_{\gcd(a_0, \ldots, a_{m-1} ,p) = 1} \left| S_{\vec{a}}(N) \right| =
O\(N^{1-\delta}\).
$$
(to see this it is enough to choose a sufficiently large $\nu$).
On the other hand, when $T$ and $N$ are close to their 
largest possible value $p^{m+1}$, that is, if
$$
T \ge N \ge p^{m+1 + o(1)}
$$
Theorem~\ref{thm:ExpSum} 
applied with $\nu =1$ gives the estimate
$$
\max_{\gcd(a_0, \ldots, a_{m-1} ,p) = 1} \left| S_{\vec{a}}(N) \right| \le 
N^{1-1/4(m+1)^2 + o(1)}.
$$
 
\subsection{Discrepancy}
\label{sec:discr}

Given a sequence $\Gamma$ of $N$ points 
\begin{equation}
\label{eq:GenSequence}
\Gamma = \left\{(\gamma_{n,0}, \ldots, \gamma_{n,m-1})_{n=0}^{N-1}\right\}
\end{equation}
in the $m$-dimensional unit cube $[0,1)^m$
it is natural to measure the level of its statistical uniformity 
in terms of the {\it discrepancy\/} $\Delta(\Gamma)$. 
More precisely, 
$$
\Delta(\Gamma) = \sup_{B \subseteq [0,1)^m}
\left|\frac{T_\Gamma(B)} {N} - |B|\right|,
$$
where $T_\Gamma(B)$ is the number of points of  $\Gamma$
inside the box
$$
B = [\alpha_1, \beta_1) \times \ldots \times [\alpha_{m}, \beta_{m})
\subseteq [0,1)^m
$$
and the supremum is taken over all such boxes, see~\cite{DrTi,KuNi}.

We recall that the discrepancy is a widely accepted 
quantitative measure 
of uniformity of distribution of sequences, and thus good pseudorandom
sequences should (after an appropriate scaling) have a small discrepancy,
see~\cite{Nied1,Nied2}.

Typically the bounds on the discrepancy of a 
sequence  are derived from bounds of exponential sums
with elements of this sequence. 
The relation is made explicit in 
 the celebrated {\it Koksma--Sz\"usz
inequality\/}, see~\cite[Theorem~1.21]{DrTi},
which we  present it in the following form.

\begin{lemma}
\label{lem:Kok-Szu} For any
integer $L > 1$ and any  sequence $\Gamma$ of $N$ points~\eqref{eq:GenSequence}
the discrepancy $\Delta(\Gamma)$
satisfies the following bound:
$$
\Delta(\Gamma)< O \( \frac{1}{L}
+ \frac{1}{N}\sum_{\substack{|a_0|, \ldots, |a_{m-1}|\le L\\
a_0^2 + \ldots + a_{m-1}^2 > 0}} \prod_{j=0}^{m-1} \frac{1}{|a_j|+1}
\left| \sum_{n=0}^{N-1} \exp \( 2 \pi i\sum_{j=0}^{m-1}a_j\gamma_{j,n} \)
\right| \).
$$
\end{lemma} 

Now, combining Lemma~\ref{lem:Kok-Szu} with the bound 
obtained in Theorem~\ref{thm:ExpSum}  
and taking $L=p-1$ we obtain:

\begin{theorem}
\label{thm:Discr} 
Let  the sequence $\{\vec{u}_n\}$ be given by~\eqref{eq:Gen}, where
the family of $m+1$ polynomials 
$\cF=\{f_0,\ldots,f_m\} \in \F_p[X_{0},\ldots,X_{m}]$ 
of total degree $d \ge 2$ is  of the form~\eqref{eq:Polys}, 
satisfying the conditions~\eqref{eq:Cond1}, \eqref{eq:Cond2} and~\eqref{eq:Cond3}, 
and such that $s_{0,1} \ldots s_{m-1,m} \ne 0$.
Assume that   $\{\vec{u}_n\}$ is purely
periodic with period $T$. Then for any fixed integer $\nu\ge 1$, 
and any positive integer $N \le T$, 
the discrepancy $D_N$ of the sequence 
$$
\(\frac{u_{n,0}}{p}, \ldots, \frac{u_{n,m-1}}{p}\),
\qquad n = 0,\ldots, N-1,
$$
satisfies the bound 
$$
D_N =
O\(p^{\alpha_{m,\nu}}N^{-\beta_{m,\nu}} (\log p)^m\)
$$
where
$$
\alpha_{m,\nu}=\frac{2m^2 + 2m\nu + 2m +\nu}{4\nu(m+\nu)} \qquad \text{and} 
\qquad \beta_{m,\nu} = \frac{1}{2\nu}
$$
and the implied constant depends  only on $d$, $m$
and $\nu$.
\end{theorem}

We remark that the same comments at the end of Section~\ref{sec:ExpSum}
also apply to Theorem~\ref{thm:Discr}. 

\section{Remarks and Open Questions}

We recall that the dynamical degree $\dd \cF$  of the 
polynomial system $\cF$ and of the associated 
affine map $\cF: \F^r \to \F^r$ 
is 
defined as
$$
\dd \cF = \lim_{k\to \infty}\(\deg \cF^{(k)}\)^{1/k}, 
$$
where $ \cF^{(k)}$ is the $k$th iteration of $\cF$
(and $\deg \cF^{(k)}$ is the largest degree of its components), 
see~\cite[Section~7.1.3]{Silv1}.
We note that the polynomial systems $\cF$ which we 
have constructed in~\eqref{eq:Polys} satisfy $\dd \cF = 1$
under the conditions~\eqref{eq:Cond1}, \eqref{eq:Cond2} and~\eqref{eq:Cond3}.
Furthermore, for any nonlinear polynomial system $\cF$  with  $\dd \cF = 1$
one can obtain an improvement of the generic bounds on
the corresponding exponential sums and the discrepancy of 
the generated sequences. However the actual improvement depends
on the speed of the convergence. 

One of the attractive choices of polynomials~\eqref{eq:Polys}, 
which leads to a very fast 
pseudorandom number generator is
$$g_i(X_{i+1},\ldots,X_m) = X_{i+1} \qquad \text{and}
\qquad h_i(X_{i+1},\ldots,X_m) = a_{i}
$$
for some constants $a_i \in \F_p$, $i = 0, \ldots, m-1$. 
The corresponding sequence of vectors is generated
at the cost of one multiplication per component. 
This naturally leads to a question of studying the 
periods of such sequences generated by such polynomial
dynamical systems.

We also note that it is natural to consider 
joint distribution of several consecutive vectors
$$
\(\vec{u_n}, \ldots, \vec{u_{n+s-1}}\), \qquad n =0,1,\ldots\,,
$$
in the $sm$-dimensional space. It seems that our method
(with some minor adjustments) 
can be  applied  to derive an  appropriate variant of 
Corollary~\ref{cor:Lin Indep} which is needed for such a result.

One of the possible ways to improve our results, 
is to construct  special polynomials $\cF = \{f_0, \ldots ,f_{r-1}\}$ 
such that linear combinations of their iterations, of the
type which appear in the proof of Theorem~\ref{thm:ExpSum}, 
satisfy the condition of the Deligne bound~\cite{Del},
that is, have a nonsingular highest form. In fact even some 
partial control over the dimension of the singularity
locus of this highest form may already lead to better estimates
via results of Katz~\cite{Katz}.

Finally, obtaining stronger results ``on average''
over all initial values $\vec{w_0} \in \F_p^{m+1}$ is an interesting and challenging 
question. It is possible that some of the arguments of~\cite{NiSh3}
may be applied to this problem.

\section*{Acknowledgement}

The authors  would  like to thank  Markus Brodmann, Joachim Rosenthal, 
Joe Silverman and Arne Winterhof for valuable comments and providing 
additional references.

During the preparation of this paper,  
A.~O. was supported in part by 
the Swiss National Science Foundation   Grant 121874  and I.~S. by
the  Australian Research Council Grant
DP0556431.


\begin{thebibliography}{99}


\bibitem{Arn1} V. I. Arnold,
`Fermat-Euler dynamical systems and the statistics of arithmetics
of geometric progressions',
{\it Func. Analysis Appl.\/}, {\bf 37} (2003), 1--15.

\bibitem{Arn2} V. I. Arnold,
`Number-theoretic turbulence in Fermat-Euler arithmetics and large Young
diagrams geometry statistics',
{\it J. Math. Fluid Mech.\/}, {\bf 7} (2005), S4--S50.


\bibitem{Arn3} V. I. Arnold,
`Ergodic and arithmetical properties of geometrical progression's
dynamics and of its orbits',    
{\it Moscow Math. J.\/}, {\bf 5} (2005), 5--22.

\bibitem{BGGS1} S. R. Blackburn, D. Gomez-Perez, J. Gutierrez and
I. E. Shparlinski, `Predicting the inversive generator',
{\it  Lect. Notes in Comp. Sci.\/}, 
Springer-Verlag, Berlin, {\bf 2898} (2003), 264--275.

\bibitem{BGGS2}  S. R. Blackburn,  D.  Gomez-Perez, J. Gutierrez and
I. E. Shparlinski, `Predicting nonlinear pseudorandom number generators', 
{\it Math. Comp.\/},  {\bf 74} (2005), 1471--1494.

\bibitem{Bour}
 J. Bourgain, `Mordell's exponential sum estimate revisited', 
{\it J. Amer. Math. Soc.\/}, {\bf 18}
(2005), 477--499.


\bibitem{Chang}
M.-C. Chang, `On a problem of Arnold on uniform distribution',
{\it J. Funcional Analysis\/},  {\bf 242} (2007),  272--280.
 
\bibitem{ChouShp} W.-S. Chou and I. E. Shparlinski,
`On the cycle structure of repeated exponentiation modulo a prime', 
{\it J. Number  Theory\/},  {\bf 107} (2004), 345--356. 

\bibitem{CohHac}  S. D. Cohen and D. Hachenberger, 
`The dynamics of linearized polynomials',
{\it Proc. Edinburgh Math. Soc.\/},
{\bf 43} (2000), 113--128. 


\bibitem{CJLS}
O. Col{\'o}n-Reyes, A. S. Jarrah,  R. Laubenbacher 
and B. Sturmfels, `Monomial dynamical systems over finite fields',
{\it Complex Systems\/},  {\bf 16}  (2006), 333--342. 

  
\bibitem{ContShp}     S. Contini and I. E. Shparlinski,
`On Stern's attack against secret
truncated linear congruential generators',
{\it Lect. Notes in Comp. Sci.\/}, Springer-Verlag,  
Berlin, {\bf 3574} (2005), 52--60.

\bibitem{Del} P. Deligne, `Applications de la formule des traces aux
sommes trigonom\'etriques', {\it Lect. Notes in Mathematics\/},
Springer-Verlag, Berlin, {\bf 569} (1977), 168--232.

\bibitem{DrTi} M. Drmota and R. Tichy, 
{\it Sequences, discrepancies and applications\/},
Springer-Verlag, Berlin, 1997.

\bibitem{EvWa} G. R. Everest  and  T. Ward,
{\it Heights of polynomials and entropy in algebraic
dynamics\/}, Springer-Verlag,  London, 1999. 
  
\bibitem{FomZel} 
S. Fomin and A. Zelevinsky, `The Laurent phenomenon',
{\it Adv. in Appl. Math.\/}, {\bf 28} (2002), 119--144.

\bibitem{FHS} J. B. Friedlander, J. Hansen and
I.~E.~Shparlinski,  `On character sums with exponential functions',
{\it Mathematika\/},  {\bf 47} (2000), 75--85.

\bibitem{FPS} J. B. Friedlander,
C.~Pomerance and I.~E.~Shparlinski, `Period of the
power generator and small values of Carmichael's
function', {\it Math. Comp.\/}, {\bf 70} (2001), 1591--1605.

\bibitem{FrSh} J. B. Friedlander and
I.~E.~Shparlinski,  `On the
distribution of the power generator',
{\it Math. Comp.\/}, {\bf 70} (2001), 1575--1589.

\bibitem{FHKLS}
A. M. Frieze, J. H{\aa}stad, R. Kannan, J. C. Lagarias and A.
Shamir, `Reconstructing truncated integer variables
satisfying linear congruences', {\it SIAM J. Comp.\/}, {\bf 17} (1988),
 262--280.
  
\bibitem{GNS} F.  Griffin, H. Niederreiter and  I.~E.~Shparlinski,
`On the  distribution
of nonlinear recursive congruential pseudorandom  numbers of higher orders',
{\it Lect. Notes in Comp. Sci.\/}, Springer-Verlag,  
Berlin, {\bf 1719} (1999), 87--93.

\bibitem{GoGuIb}   D.  Gomez-Perez, J. Gutierrez and
{\'A}. Ibeas,  `Attacking the Pollard generator',
{\it  IEEE Trans.
Inform. Theory\/}, {\bf 52} (2006), 5518--5523.

\bibitem{GG} J. Gutierrez and D. Gomez-Perez,
`Iterations of multivariate polynomials and discrepancy of pseudorandom numbers',
{\it Lect. Notes in Comp. Sci.\/}, Springer-Verlag,  
Berlin, {\bf 2227} (2001),  192--199.

\bibitem{GoGuSh}   D.  Gomez-Perez, J. Gutierrez 
and I. E, Shparlinski,
`Exponential sums with Dickson polynomials', 
{\it Finite Fields Appl.\/}, {\bf 12} (2006),  16--25.

\bibitem{GuIb}   J. Gutierrez and
{\'A}. Ibeas,  `Inferring sequences produced by a linear congruential
generator on elliptic curves missing high-order bits',
{\it Designs, Codes and Cryptography\/}, 
{\bf 41} (2007), 199--212.

\bibitem{GuWin} J. Gutierrez and A. Winterhof,  
`Exponential sums of nonlinear congruential pseudorandom number generators 
with redei functions', 
{\it Finite Fields Appl.\/}, {\bf 14} (2008), 410--416.

\bibitem{Jon} R. Jones, `The density of prime divisors in 
the arithmetic dynamics of quadratic polynomials',
{\it J. Lond. Math. Soc.\/}, 
{\bf  78} (2008), 523--544.

\bibitem{JS} A. Joux and J. Stern, `Lattice reduction: 
A toolbox for the cryptanalyst',
{\it J. Cryptology\/}, {\bf 11} (1998), 161--185.

\bibitem{Katz} N. Katz, `Estimates for ``singular'' exponential sums',
{\it International Mathematics Research Notices\/}, 
\textbf{16} (1999), 875--899.

\bibitem{Kraw} H. Krawczyk, `How to predict congruential generators',
{\it  J. Algorithms\/}, {\bf 13} (1992), 527--545.

\bibitem{KuNi}
L.~Kuipers and H.~Niederreiter, {\it Uniform distribution of
sequences\/}, Wiley-Intersci., New York-London-Sydney, 1974.

\bibitem{Lag} J. C. Lagarias, `Pseudorandom number 
generators in cryptography and number
theory', {\it Proc. Symp. in Appl. Math.\/}, 
Amer. Math. Soc., Providence, RI,
{\bf 42} (1990), 115--143.

\bibitem{LN} R. Lidl and H. Niederreiter, {\it Finite fields\/},
Cambridge University Press, Cambridge, 1997.

\bibitem{Marc} S. Marcello, `Sur la dynamique 
arithm\'etique des automorphismes de l'espace
affine', {\it Bull. Soc. Math. France\/}, {\bf 131} 
(2003), 229--257.  

\bibitem{MartPom} G.~Martin  and C.~Pomerance, 
`The iterated Carmichael $\lambda$-function and the 
number of cycles of the power generator', 
{\it Acta Arith.\/}, {\bf 118} (2005), 305--335.


\bibitem{Nied1} H. Niederreiter, `Quasi-Monte Carlo methods and pseudo-random
numbers', {\it Bull. Amer. Math. Soc.\/}, {\bf 84} (1978), 957--1041.

\bibitem{Nied2}
H. Niederreiter, {\it Random number
generation and Quasi--Monte Carlo methods\/},  SIAM Press, 1992.


\bibitem{NiSh1} H. Niederreiter and I. E. Shparlinski,
`On the distribution and lattice structure of nonlinear congruential pseudorandom
numbers', {\it Finite Fields and Their Appl.\/}, {\bf 5} (1999), 246--253.

\bibitem{NiSh2}
H. Niederreiter and I. E. Shparlinski,
`On the distribution of inversive congruential pseudorandom
numbers in parts of the period',
{\it Math. Comp.\/}, {\bf 70} (2001), 1569--1574.
 

\bibitem{NiSh3}
H. Niederreiter and I. E. Shparlinski,
`On the average distribution of inversive
pseudorandom numbers',
{\it Finite Fields and Their Appl.\/},  
{\bf 8} (2002), 491--503.

\bibitem{NiSh4} H. Niederreiter and I. E. Shparlinski,
`Dynamical systems generated by rational functions',
{\it Lect. Notes in Comp. Sci.\/}, Springer-Verlag,  
Berlin, {\bf 2643} (2003), 6--17.




\bibitem{NiWi}
H. Niederreiter and A. Winterhof,  
`Exponential sums for nonlinear recurring sequences', 
{\it Finite Fields Appl.\/}, {\bf 14} (2008),   59--64. 

\bibitem{Shp} I.~E.~Shparlinski, `On some dynamical systems in  finite fields 
and residue rings',
{\it Discr. and Cont. Dynam. Syst., Ser.A\/},  {\bf 17} (2007), 901--917.

\bibitem{Silv1}
J. H. Silverman, {\it The arithmetic of dynamical systems\/},
Springer, New York, 2007.

\bibitem{Silv2}
J. H. Silverman, `Variation of periods modulo $p$ in arithmetic dynamics',
{\it New York J. Math.\/},  {\bf 14}  (2008), 601--616.  

\bibitem{TopWin} A. Topuzo{\v g}lu and
A. Winterhof, `Pseudorandom  sequences', {\it Topics in Geometry,
Coding Theory and Cryptography\/}, Springer-Verlag, 2006, 135--166.

\bibitem{VaSha} T. Vasiga and J. O.  Shallit, 
`On the iteration of certain quadratic maps over ${\mathrm{GF}(p)}$', 
{\it  Discr. Math.\/},  {\bf 277} (2004), 219--240. 

\end{thebibliography}
\end{document}